\documentclass[10pt]{amsart}

\usepackage[colorlinks]{hyperref}
\usepackage{color,graphicx,shortvrb}
\usepackage[latin 1]{inputenc}
\usepackage[dvipsnames]{xcolor}

\usepackage[active]{srcltx} 

\usepackage{enumerate}
\usepackage{amssymb}
\usepackage{tikz-cd}

\usepackage [ all ]{xy}

\newtheorem{theorem}{Theorem}[section]

\newtheorem{corollary}[theorem]{Corollary}

\newtheorem{definition}[theorem]{Definition}
\newtheorem{lemma}[theorem]{Lemma}

\newtheorem{proposition}[theorem]{Proposition}

\def\J#1#2#3{ \left\{ #1,#2,#3 \right\} }

\def\11{\textbf{$1$}}
\def\11b#1{\mathbf{1}_{_{#1}}}
\def\CC{{\mathbb{C}}}

\DeclareGraphicsExtensions{.jpg,.pdf,.png,.eps}

\begin{document}

\title[Maps preserving triple transition pseudo-probabilities]{Maps preserving triple transition pseudo-probabilities}

\author[A.M. Peralta]{Antonio M. Peralta}

\address{Instituto de Matemáticas de la Universidad de Granada (IMAG). Departamento de An{\'a}lisis Matem{\'a}tico, Facultad de	Ciencias, Universidad de Granada, 18071 Granada, Spain.}

\email{aperalta@ugr.es}


\subjclass[2010]{Primary 47B49, 46L60, 47N50  Secondary  81R15, 17C65}

\keywords{Wigner theorem, minimal partial isometries, minimal tripotents, socle, triple transition psedudo-probability, preservers, Cartan factors, spin factors, triple isomorphism}

\date{}

\begin{abstract} Let $e$ and $v$ be minimal tripotents in a JBW$^*$-triple $M$.  We introduce the notion of triple transition pseudo-probability from $e$ to $v$ as the complex number $TTP(e,v)= \varphi_v(e),$ where $\varphi_v$ is the unique extreme point of the closed unit ball of $M_*$ at which $v$ attains its norm. In the case of two minimal projections in a von Neumann algebra, this correspond to the usual transition probability. We prove that every bijective transformation $\Phi$ preserving triple transition pseudo-probabilities between the lattices of tripotents of two atomic JBW$^*$-triples $M$ and $N$ admits an extension to a bijective {\rm(}complex{\rm)} linear mapping between the socles of these JBW$^*$-triples. If we additionally assume that $\Phi$ preserves orthogonality, then $\Phi$ can be extended to a surjective (complex-)linear {\rm(}isometric{\rm)} triple isomorphism from $M$ onto $N$. In case that $M$ and $N$ are two spin factors or two type 1 Cartan factors we show, via techniques and results on preservers, that every bijection preserving triple transition pseudo-probabilities between the lattices of tripotents of $M$ and $N$ automatically preserves orthogonality, and hence admits an extension to a triple isomorphism from $M$ onto $N$.  
\end{abstract}

\maketitle
\thispagestyle{empty}

\section{Introduction}

The available mathematical models for quantum mechanic make use of complex Hilbert spaces to define the states of a quantum system. Given a complex Hilbert space $H$, the normal state space of $S(H)$ is identified, via trace duality, with those positive norm-one elements (states) in the predual of the von Neumann algebra, $B(H),$ of all bounded linear operators on $H$. Each observable is associated with a self-adjoint operator in $A\in B(H),$ and its expected value on the the system in state $p$ is $A(p)=\hbox{tr}(Ap),$ where $tr(.)$ stands for the usual trace on $B(H)$. The elements in $S(H)$ are called the normal states of a quantum system associated to the Hilbert space $H$. The extreme points of $S(H),$ as a convex set inside the closed unit ball of $B(H)_*,$ are called pure states, and they can be also identified with rank-one projections on $H$. The set of all rank-one projections on $H$ will be  denoted by $\mathcal{P}_1 (H)$, while $\mathcal{P}(H)$ or $\mathcal{P} (B(H))$ will stand for the set of all (orthogonal) projections on $H$.\smallskip

If two pure states are represented by the minimal projections $p= \xi\otimes \xi$ and $q= \eta\otimes \eta,$ with $\xi$ and $\eta$ in the unit sphere of $H$, according to Born's rule, the \emph{transition probability} from $p$ to $q$ is defined as
$$TP(p,q)=\hbox{tr}(p q) = \hbox{tr}(pq^*) = \hbox{tr}(q p^*) = |\langle \xi, \eta\rangle|^2.$$ Here and along this note, for $\xi$ in another complex Hilbert space $K$ and $\eta\in H$, the symbol $\xi\otimes \eta$ will stand for the operator from $H$ to $K$ defined by $\xi\otimes \eta (\zeta) := \langle \zeta, \eta\rangle \xi.$\smallskip

A bijective map $\Phi : \mathcal{P}_1 (H) \to \mathcal{P}_1 (H)$ is called a symmetry transformation or a Wigner symmetry if it preserves the transition probability between minimal projections, that is, $$TP (\Phi(p),\Phi(q)) = \hbox{tr} ( \Phi(p) \Phi(q)) = \hbox{tr} (pq) = TP(p,q),\hbox{ for all } (p,q \in \mathcal{P}_1 (H)).$$ 

A linear (respectively, conjugate-linear) mapping $u : H\to H$ is called a unitary (respectively, an anti-unitary) if $u u ^* = u^* u =1$. The celebrated Wigner's theorem admits the following statement:

\begin{theorem}\label{t Wigner}{\rm(Wigner theorem, \cite{Wig31}, \cite[page 12]{Molnar85})} Let $H$ be a complex Hilbert space. A bijective mapping $\Phi : \mathcal{P}_1(H) \to \mathcal{P}_1(H)$ is a symmetry
transformation if and only if there is an either unitary or anti-unitary operator $u$ on $H$, unique up to multiplication by a unitary scalar, such that $\Phi (p) = u p u^*$ for all $p\in \mathcal{P}_1(H)$. Furthermore, the real linear {\rm(}actually complex-linear or conjugate-linear{\rm)} mapping $T: B(H) \to B(H),$ $T(x) = u x u^*$ is a $^*$-automorphism whose restriction to $ \mathcal{P}_1(H)$ coincides with $\Phi$.
\end{theorem}

It is known (see, for example, \cite[\S 4 and 6]{CasdeVilahtiLevrero97}) that for a complex Hilbert space $H$ with dim$(H)\geq 3$, the following mathematical models employed in the Hilbert space formulation of quantum mechanics are equivalent: \begin{enumerate}[$(M.1)$]
	\item The set $\mathbf{P}$ of pure states on $H$ (which algebraically corresponds to the set $\mathcal{P}_1(H)$) whose automorphisms are the bijections preserving transition probabilities.
	\item The orthomodular lattice $\mathbf{L}$ of closed subspaces of $H$, or equivalently, the lattice of all projection in $B(H),$ where the automorphisms are the bijections preserving orthogonality and order. 
\end{enumerate} 

The equivalence of these two models implies that if dim$(H)
\geq 3$, every bijection $\Phi: \mathcal{P} (B(H))\to  \mathcal{P} (B(H))$ preserving the partial ordering and orthogonality in both directions is given by a real linear $^*$-automorphism on $B(H)$ determined either by a unitary or by an anti-unitary operator on $H$ (cf. \cite[\S 2.3 and Proposition 4.9]{CasdeVilahtiLevrero97}).\smallskip

The lattice of projections in $B(H)$ is a subset of the strictly bigger lattice of partial isometries in $B(H)$. We recall that an element $e$ in $B(H)$ is a partial isometry if $e e^*$ (equivalently, $e^* e$) is a projection. Partial isometries are also called tripotents since an element $e$ is a partial isometry if and only if $e e^* e =e$. Let the symbol $\mathcal{PI}(H)= \mathcal{U}(B(H))$ stand for the set of all partial isometries on $H$. We shall write $\mathcal{PI}_1(H) = \mathcal{U}_{min} (B(H))$ for the set of all rank--1 or minimal partial isometries on $H$. We say that $e,v\in \mathcal{U}(B(H))$ are orthogonal if and only if $\{ ee^*,vv^*\}$ and $\{ e^*e,v^* v\}$ are two sets of orthogonal projections. The standard partial ordering on $\mathcal{U}(B(H))$ is defined in the following terms: $ e\leq u$ if $u-e$ is a partial isometry orthogonal to $e$.\smallskip

L. Moln{\'a}r seems to be the first author in considering a Wigner type theorem for bijections on the lattice of partial isometries of $B(H)$ preserving the partial order and orthogonality in both directions. 

\begin{theorem}\label{t Molnar 2002}\cite[Theorem 1]{Molnar2002} Let $H$ be a complex Hilbert space with dim$(H)\geq 3$. Suppose that $\Phi : \mathcal{U}(B(H))\to \mathcal{U}(B(H))$ is a bijective transformation which preserves the partial ordering and the orthogonality between partial isometries in both directions. If $\Phi$ is continuous {\rm(}in the operator norm{\rm)} at a single element of $\mathcal{U}(B(H))$ different from $0$, then $\Phi$ extends to a real-linear triple isomorphism.
\end{theorem}

During the mini-symposium \emph{``Research on preserver problems on Banach algebras and related topics''} held at RIMS (Research Institute for Mathematical Sciences), Kyoto University on October 25--27, 2021, the author of this note presented the following generalization of the previous theorem to the case of atomic JBW$^*$-triples (i.e. JB$^*$-triples which are $\ell_{\infty}$-sums of Cartan factors).  

\begin{theorem}\label{t order and orthogonality preservers atomic JBWtriples}{\rm\cite[Theorem 6.1]{FriPe2021}} Let $\displaystyle M= \bigoplus_{i\in I}^{\ell_{\infty}} C_i$ and  $\displaystyle N = \bigoplus_{j\in J}^{\ell_{\infty}} \tilde{C}_j$ be atomic JBW$^*$-triples, where $C_i$ and $C_j$ are Cartan factors with rank $\geq 2$. Suppose that $\Phi : \mathcal{U}(M) \to \mathcal{U}(N)$ is a bijective transformation which preserves the partial ordering in both directions and orthogonality between tripotents. We shall additionally assume that $\Phi$ is continuous at a tripotent $u = (u_i)_i$ in $M$ with $u_i\neq 0$ for all $i$ {\rm(}or we shall simply assume that $\Phi|_{\mathbb{T} u}$ is continuous at a tripotent $(u_i)_i$ in $M$ with $u_i\neq 0$ for all $i${\rm)}. Then there exists a real linear triple isomorphism $T: M\to N$ such that $T(w) = \Phi(w)$ for all $w\in \mathcal{U} (M)$. 
\end{theorem}

It should be remarked that the hypothesis concerning the ranks of the Cartan factors in the previous theorem cannot be relaxed (cf. \cite[Remark 3.6]{FriPe2021}). Anyway, the validity of the result for rank-2 Cartan factors is undoubtedly an advantage.\smallskip

Back to the essence of Wigner theorem expressed in Theorem \ref{t Wigner}, we find the following contribution by L. Moln{\'a}r.

\begin{theorem}\label{t Molnar minimal pi}\cite[Theorem 2]{Molnar2002} Let $\Phi: \mathcal{U}_{min} (B(H))\to \mathcal{U}_{min} (B(H))$ be a bijective mapping satisfying \begin{equation}\label{equation in Molnar's teorem} \hbox{tr} ({\Phi(e)}^* \Phi(v)) = \hbox{tr}(e^* v), \hbox{ for all } e,v\in \mathcal{U}_{min} (B(H)).
	\end{equation}
Then $\Phi$ extends to a surjective complex-linear isometry. Moreover, one of the following statements holds: \begin{enumerate}[$(a)$]
		\item there exist unitaries $u, w$ on $H$ such that $\Phi(e) = u e w$ {\rm(}$e \in \mathcal{U}_{min}(B(H))${\rm)};
		\item there exist anti-unitaries $u, w$ on $H$ such that $\Phi(e) = u e^* w$ {\rm(}$e \in \mathcal{U}_{min}(B(H))${\rm)}.
	\end{enumerate}
\end{theorem}

Let us observe that for each minimal partial isometry $e$ in $B(H),$ the functional $\varphi_e (x) = \hbox{tr} (e^* x)$ is the unique extreme point of the closed unit ball of $B(H)_*$, the predual of $B(H)$, at which $e$ attains its norm. A similar property holds in the wider setting of JBW$^*$-triples (see subsection \ref{subsect: definitions} for details and definitions). Namely, for each minimal tripotent $e$ in a JBW$^*$-triple, $M,$ there exists a unique pure atom (i.e. an extreme point of the closed unit ball of $M_*$) $\varphi_e$ at which $e$ attains its norm and the corresponding Peirce-2 projection writes in the form $P_2 (e) (x) = \varphi_e(x) e$ for all $x\in M$ (cf. \cite[Proposition 4]{FriRu85}). The mapping $$\mathcal{U}_{min} (M)\to \partial_{e} (\mathcal{B}_{M_*}), \ \  e\mapsto \varphi_e $$ is a bijection from the set of minimal tripotents in $M$ onto the set of pure atoms of $M$. Given two minimal tripotents $e$ and $v$ in a JBW$^*$-triple $M$, we define the \emph{triple transition pseudo-probability} from $e$ to $v$ as the complex number given by $TTP(e,v)= \varphi_v(e).$ So, the hypothesis \eqref{equation in Molnar's teorem} in Theorem \ref{t Molnar minimal pi} is equivalent to say that $\Phi$ preserves triple transition pseudo-probabilities. In the case of $B(H)$, the triple transition pseudo-probability between two minimal projections is precisely the usual transition probability. We shall show that this pseudo-probability is symmetric in the sense that $TTP(e,v)= 
\overline{TTP(v,e)},$ for every couple of minimal tripotents $e,v\in M$. \smallskip

We shall see below that the triple transition pseudo-probability between any two minimal projections $p$ and $q$ in a von Neumann algebra $W$ is zero if and only if $p$ and $q$ are orthogonal (i.e. $p q =0$). The same equivalence does not necessarily hold when projections are replaced with tripotents or partial isometries, for example, the partial isometries $e = \left(                                                                             \begin{array}{cc} 1 & 0 \\
	0 & 0 \\
\end{array}
\right)$ and $v = \left(
\begin{array}{cc}
	0 & 1 \\
	0 & 0 \\
\end{array}
\right)$ are not orthogonal in $M_2(\mathbb{C})$, but $TTP(e,v)=0$. This is a theoretical handicap for the triple transition pseudo-probability. However, despite Theorem 
\ref{t order and orthogonality preservers atomic JBWtriples} does not hold for rank-one JB$^*$-triples (cf. \cite[Remark  3.6]{FriPe2021}), every (non-necessarily surjective) mapping between the lattices of tripotents of two rank-one JB$^*$-triples preserving triple transition pseudo-probabilities  always admits an extension to a linear and isometric triple homomorphism between the corresponding JB$^*$-triples (see Proposition \ref{p rank-one Cartan factors}). This will be obtained by an application of a theorem of Ding on the extension of isometries on the unit sphere of a Hilbert space \cite{Ding2002}.\smallskip

In Theorem \ref{t bijections preserving triple transition pseudo-probabilities} we establish that if $M$ and $N$ are atomic JBW$^*$-triples and $\Phi : \mathcal{U}_{min}(M) \to \mathcal{U}_{min}(N)$ is a bijective transformation preserving triple transition pseudo-probabilities between the sets of minimal tripotents, then there exists a bijective {\rm(}complex{\rm)} linear mapping $T_0$ from the socle of $M$ onto the socle of $N$ whose restriction to $\mathcal{U}_{min} (M)$ is $\Phi$, where the socle of a JB$^*$-triple is the subspace linearly generated by its minimal tripotents. If we additionally assume that $\Phi$ preserves orthogonality, then we prove the existence of a surjective (complex-)linear {\rm(}isometric{\rm)} triple isomorphism from $M$ onto $N$ extending the mapping $\Phi$ (cf. Corollary \ref{c bijections preserving triple transition pseudo-probabilities and orthogonality}). \smallskip

Due to the just commented result, the natural question is whether every bijection between the sets of minimal tripotents in two atomic JBW$^*$-triples preserving triple transition pseudo-probabilities must automatically preserve orthogonality among them. The rest of the paper is devoted to present a couple of positive answers to this problem in the case of spin and type 1 Cartan factors. \smallskip

Section \ref{sec: spin factors} is devoted to study bijections preserving triple transition pseudo-proba-bilities between the sets of minimal tripotents in two spin factors. We shall show that any such bijection preserves orthogonality, and hence admits an extension to a triple isomorphism between the spin factors (see Theorem \ref{t bijections preserving triple transition pseudo-probabilities spin factors})). The proof is based on an remarkable results on preservers, due to J. Chmieli{\'n}ski,  asserting that a non-vanishing mapping between two inner product spaces is linear and preserves orthogonality in the Euclidean sense if and only if it is a positive multiple of a linear isometry \cite[Theorem 1]{Chmielinski05}.\smallskip

In section \ref{sec: type 1 Cartan factors} we also establish a positive answer to the problem stated above in the case of a bijection between the sets of minimal tripotents in two type 1 Cartan factors (see Theorem \ref{t bijections preserving triple transition pseudo-probabilities type 1 Cartan factors}). On this occasion, our arguments run closer to those given by Moln{\'a}r in the proof of Theorem \ref{t Molnar minimal pi}. For this purpose we shall establish a variant of several results previously explored by M. Marcus, B.N. Moyls \cite{MarcusMoyls59}, R. Westwick \cite{Westwick67} and M. Omladi\v{c} and P. \v{S}emrl \cite{OmladicSemrl93}. We concretely prove, in Theorem \ref{t extension of OmladicSemrl}, that for each linear bijection $\Phi : soc(B(H_1,K_1))\to soc(B(H_2,K_2))$ preserving rank-one operators in both directions, where $H_1,H_2, K_1$ and $K_2$ are complex Hilbert spaces with dimensions $\geq 2$, one of the next statements holds: \begin{enumerate}[$(a)$]\item either there are bijective linear mappings $u:K_1\to K_2,$ and $v:H_1\to H_2$ such that $\Phi (\xi \otimes \eta) = u (\xi\otimes \eta) v= u(\xi) \otimes v(\eta)$ {\rm(}$\xi\in K_1,\eta \in H_1${\rm)};
	\item or there are bijective conjugate-linear mappings $u:H_1\to K_2,$ $v:K_1\to H_2$ such that $\Phi (\xi \otimes \eta)= u (\xi\otimes \eta)^* v $ $= u (\eta\otimes \xi) v= u(\eta) \otimes v(\xi)$ {\rm(}$\xi\in K_1,$ $\eta \in H_1${\rm)}.
\end{enumerate}  

Let us finish this introduction with a kind of announcement or statement of intentions, it would be desirable to find a positive argument to prove that every bijection between the sets of minimal tripotents in two atomic JBW$^*$-triple automatically preserves orthogonality. Perhaps a more general point of view could provide a better understanding. At the present moment it seems a open problem. Some other additional questions also arise after this first study on triple transition pseudo-probabilities.

\subsection{Definitions and terminology}\label{subsect: definitions}\ \smallskip

The model which motivated the study of C$^*$-algebras is the space $B(H),$ of all bounded linear operators on a complex Hilbert space $H$. Left and right weak$^*$ closed ideals of $B(H)$ are precisely subspaces of the form $ B(H)p$ and $pB(H),$ respectively, where $p$ is a projection in $B(H)$. These ideals are identified with subspaces of operators of the form $B(p(H), H)$ and $B(H, p(H))$. However, given two complex Hilbert spaces $H$ and $K$ (where we can always assume that $K$ is a closed subspace of $H$), the Banach space $B(H,K)$, of all bounded linear operators from $H$ to $K$, is not, in general, a C$^*$-subalgebra of some $B(H)$. Despite of this handicap, $B(H,K)$ is stable under products of the form \begin{equation}\label{eq triple product JCstar triple} \{x,y,z\} = \frac12\left( x y^* z + z y^* x\right)\ \ (x,y,z\in B(H,K)).
\end{equation}  Closed complex-linear subspaces of $B(H,K)$ which are closed for the triple product defined in \eqref{eq triple product JCstar triple} were called \emph{J$^*$-algebras} by L. Harris in \cite{Harris74,Harris81}. J$^*$-algebras include, in particular, all C$^*$-algebras, all JC$^*$-algebras, all complex Hilbert spaces, and all ternary algebras of operators. Harris also proved that the open unit ball of every J$^*$-algebra enjoys the interesting holomorphic property of being a bounded symmetric domain (see \cite[Corollary 2]{Harris74}). In \cite{BraKaUp78}, R. Braun, W. Kaup and H. Upmeier extended Harris' result by showing that the open unit ball of every (unital) JB$^*$-algebra satisfies the same property.\smallskip

If the holomorphic-property "being a bounded symmetric domain" is employed to classify the open unit balls of complex Banach spaces, the definitive result is due to W. Kaup, who in his own words \emph{``introduced the concept of a JB$^*$-triple and showed that every bounded symmetric domain in a complex Banach space is biholomorphically equivalent to the open unit ball of a JB$^*$-triple, and in this way, the category of all bounded symmetric domains with base point is equivalent to the category of JB$^*$-triples}'' (see \cite{Ka83}).\smallskip

A complex Banach space $E$ is called a \emph{JB$^*$-triple} if it admits a continuous triple product $\J \cdot\cdot\cdot :
E\times E\times E \to E,$ which is symmetric and bilinear in the first and third variables, conjugate-linear in the middle one,
and satisfies the following axioms:
\begin{enumerate}[{\rm (a)}] \item (Jordan identity)
$$L(a,b) L(x,y) = L(x,y) L(a,b) + L(L(a,b)x,y)
 - L(x,L(b,a)y)$$ for $a,b,x,y$ in $E$, where $L(a,b)$ is the operator on $E$ given by $x \mapsto \J abx;$
\item $L(a,a)$ is a hermitian operator with non-negative spectrum for all $a\in E$;
\item $\|\{a,a,a\}\| = \|a\|^3$ for every $a\in E$.\end{enumerate}

The first examples of JB$^*$-triples include C$^*$-algebras and $B(H,K)$ spaces with respect to the triple product given in \eqref{eq triple product JCstar triple}, the latter are known as \emph{Cartan factors of type 1}.\smallskip

There are six different types of Cartan factors, the first one has been introduced in the previous paragraph. In order to define the next two types, let $j$ be a conjugation (i.e. a conjugate-linear isometry or period 2) on a complex Hilbert space $H$. We consider a linear involution on $B(H)$ defined by $x\mapsto x^t:=jx^*j$. \emph{Cartan factors of type 2 and 3} are the JB$^*$-subtriples of $B(H)$ of all $t$-skew-symmetric and $t$-symmetric operators, respectively.\smallskip

A \emph{Cartan factor of type 4}, also called a \emph{spin factor},\label{def spin factor} is a complex Hilbert space $M$ provided with a conjugation $x\mapsto \overline{x},$ where the triple product and the norm are defined by \begin{equation}\label{eq spin product}
\{x, y, z\} = \langle x, y\rangle z + \langle z, y\rangle  x -\langle x, \overline{z}\rangle \overline{y},
\end{equation} and \begin{equation}\label{eq spin norm} \|x\|^2 = \langle x, x\rangle  + \sqrt{\langle x, x\rangle ^2 -|
\langle x, \overline{x}\rangle  |^2},
 \end{equation} respectively (cf. \cite[Chapter 3]{Fri2005}). The \emph{Cartan factors of types 5 and 6} (also called \emph{exceptional} Cartan factors) are spaces of matrices over the eight dimensional complex algebra of Cayley numbers; the type 6 consists of all $3\times 3$ self-adjoint matrices and has a natural Jordan algebra structure, and the type 5 is the subtriple consisting of all $1\times 2$ matrices (see \cite{Ka97, Harris74, HervIs92} and the recent references \cite[\S 6.3 and 6.4]{HamKalPe20}, \cite[\S 3]{HamKalPe22determinants} for more details).\smallskip

An element $e$ in a JB$^*$-triple $E$ is called a \emph{tripotent} if $\{e,e,e\}= e$. When a C$^*$-algebra is regarded as a JB$^*$-triple with the triple product in \eqref{eq triple product JCstar triple}, tripotents and partial isometries correspond to the same elements. If we fix a tripotent $e$ in $E$, we can find a decomposition of the space in terms of the eigenspaces of the operator $L(e,e)$ which is expressed as follows:
\begin{equation}\label{Peirce decomp} {E} = {E}_{0} (e) \oplus  {E}_{1} (e) \oplus {E}_{2} (e),\end{equation} where ${
E}_{k} (e) := \{ x\in {E} : L(e,e)x = {\frac k 2} x\}$ is a subtriple of ${E}$ called the \emph{Peirce-$k$ subspace} ($k=0,1,2$). \emph{Peirce-$k$ projection} is the name given to the natural projection of ${E}$ onto ${E}_{k} (e),$ and it is usually denoted by $P_{k} (e)$. 
 The Peirce-$2$ subspace ${E}_{2} (e)$ is a unital JB$^*$-algebra with respect to the product and involution given by $x \circ_e y = \J xey$ and $x^{*_e} = \J exe,$ respectively. 
\smallskip

A tripotent $e$ in $E$ is called \emph{algebraically minimal} (respectively, \emph{complete} or \emph{algebraically maximal}) if  $E_2(e)=\CC e \neq \{0\}$ (respectively, $E_0 (e) =\{0\}$). We shall say that $e$ is a \emph{unitary tripotent} if $E_2(e) =E$. The symbols $\mathcal{U} (E)$, $\mathcal{U}_{min} (E)$, and $\mathcal{U}_{max} (E)$ will stand for the sets of all tripotents, minimal tripotents, and complete tripotents in $E$, respectively.\smallskip

A JB$^*$-triple might contain no non-trivial tripotents, that is the case of the JB$^*$-triple $C_0[0,1]$ of all complex-valued continuous functions on $[0,1]$ vanishing at $0$. However, in a JB$^*$-triple $E$ the extreme points of its closed unit ball are precisely the complete tripotents in $E$ (cf. \cite[Lemma 4.1]{BraKaUp78}, \cite[Proposition 3.5]{KaUp77} or \cite[Corollary 4.8]{EdRutt88}). Thus, every JB$^*$-triple which is also a dual Banach space contains an abundant set of tripotents. JB$^*$-triples which are additionally dual Banach spaces are called \emph{JBW$^*$-triples}. Each JBW$^*$-triple admits a unique (isometric) predual and its triple product is separately weak$^*$ continuous (cf. \cite{BarTi}).\smallskip

A JBW$^*$-triple is called \emph{atomic} if it coincides with the w$^*$-closure of the linear span of its minimal tripotents. A very natural example is given by $B(H),$ where each minimal tripotent is of the form $\xi\otimes \eta$ with $\xi,\eta$ in the unit sphere of $H$. Every Cartan factor is an atomic JBW$^*$-triple. Cartan factors are enough to exhaust all possible cases since every atomic JBW$^*$-triple is an $\ell_{\infty}$-sum of Cartan factors (cf. \cite[Proposition 2 and Theorem E]{FriRu86}).\smallskip

The notion of orthogonality between tripotents is an important concept in the theory of JB$^*$-triples. Suppose $e$ and $v$ are two tripotents in a JB$^*$-triple $E$. According to the standard notation (see, for example \cite{loos1977bounded,Batt91}) we say that $e$ is \emph{orthogonal} to $u$ ($e\perp u$ in short) if $\{e,e,u\}=0$. It is known that $e\perp u$ if and only if $\{u,u,e\}=0$ (and the latter is equivalent to any of the next statements $L(e,u) = 0,$ $L(u,e) = 0,$  $e \in E_0(u),$ $u\in E_0(e)$ cf. \cite[Lemma 3.9]{loos1977bounded}). It is worth to remark that two projections $p$ and $q$ in a C$^*$-algebra $A$, regarded as a JB$^*$-triple, are orthogonal if and only if $p q =0$  (that is, they are orthogonal in the usual sense).\smallskip

We can also speak about orthogonality for pairs of general elements in a JB$^*$-triple $E$. We shall say that $x$ and $y$ in $E$ are \emph{orthogonal} ($x\perp y$ in short) if $L(x,y)=0$ (equivalently $L(y,x)=0$, compare \cite[Lemma 1.1]{BurFerGarMarPe08} for several reformulations). Any two orthogonal elements $a$ and $b$ in JB$^*$-triple $E$ are $M$-orthogonal in a strict geometric sense, that is, $\|a+b\| = \max\{\|a\|, \|b\|\}$ (see \cite[Lemma 1.3$(a)$]{FriRu85}). \smallskip

Building upon the relation ``being orthogonal'' we can define a canonical order ``$\leq$'' on tripotents in $E$ given by $e\leq u$ if and only if $u-e$ is a tripotent and $u-e \perp e$. This partial ordering is precisely the order consider by L. Moln{\'a}r in Theorem \ref{t Molnar 2002}, and it provides an important tool in JB$^*$-triples (see, for example, the recent papers \cite{HamKalPePfi20,HamKalPePfi20Groth,HamKalPe20, HamKalPe22determinants, Ham21, HamKalPeOrder} where it plays an important role). The partial order in $\mathcal{U} (E)$ enjoys several interesting properties; for example, $e\leq u$ if and only if $e$ is a projection in the JB$^*$-algebras $E_2(e)$ (cf. \cite[Lemma 3.2]{Batt91} or \cite[Corollary 1.7]{FriRu85} or \cite[Proposition 2.4]{HamKalPe20}). In particular, if $e$ and $p$ are tripotents (i.e. partial isometries) in a C$^*$-algebra $A$ regarded as a JB$^*$-triple with the triple product in \eqref{eq triple product JCstar triple} and $p$ is a projection, the condition $e\leq p$ implies that $e$ is a projection in $A$ with $e\leq p$ in the usual order on projections (i.e. $p e = e$).\smallskip

A non-zero tripotent $e$ in $E$ is called (\emph{order}) \emph{minimal} (respectively, (\emph{order}) \emph{maximal}) if $0\neq u \leq e$ for a tripotent $u$ in $E$ implies that $u =e$ (respectively, $e \leq u$ for a tripotent $u$ in $E$ implies that $u =e$). Clearly, every algebraically minimal tripotent is (order) minimal but the reciprocal implication does not necessarily hold, for example, the unit element in $C[0,1]$ is order minimal but not algebraically minimal. In the C$^*$-algebra $C_0[0,1],$ of all continuous functions on $[0,1]$ vanishing at 0, the zero tripotent is order maximal but it is not algebraically maximal. In the setting of JBW$^*$-triples these pathologies do not happen, that is, in a JBW$^*$-triple order and algebraic maximal (respectively, minimal) tripotents coincide (cf. \cite[Corollary 4.8]{EdRutt88} and \cite[Lemma 4.7]{Batt91}).\smallskip

A triple homomorphism between JB$^*$-triples $E$ and $F$ is a linear map $T:E\to F$ such that $T\{a,b,c\} = \{T(a),T(b), T(c)\}$ for all $a,b,c\in E$. Every triple homomorphism between JB$^*$-triples is continuous \cite[Lemma 1]{BarDanHor88}. A triple isomorphism is a bijective triple homomorphism. Clearly, the inclusion $T(\mathcal{U} (E)) \subseteq \mathcal{U} (F)$ holds for each triple homomorphism $T$, while the equality  $T(\mathcal{U} (E)) =\mathcal{U} (F)$ is true for every triple isomorphism $T$. Every injective triple homomorphism is an isometry (see \cite[Lemma 1]{BarDanHor88}). Actually a deep result in the theory of JB$^*$-triples, established by W. Kaup in \cite[Proposition 5.5]{Ka83}, proves that a linear bijection between JB$^*$-triples is a triple isomorphism if and only if it is an isometry. Therefore, each triple isomorphism $T: E\to F$ induces a surjective isometry $T|_{\mathcal{U} (E)}: \mathcal{U} (E)\to \mathcal{U} (F)$ which preserves orthogonality and partial order in both directions. Similar arguments prove that the mappings $T|_{\mathcal{U}_{min} (E)}: \mathcal{U}_{min} (E)\to \mathcal{U}_{min} (F)$ and $T|_{\mathcal{U}_{max} (E)}: \mathcal{U}_{max} (E)\to \mathcal{U}_{max} (F)$ are surjective isometries.\smallskip

Along this note, the unit sphere of each normed space $X$ will be denoted by $S_{_X}$, and we shall write $\mathbb{T}$ for $S_{_\mathbb{C}}$.

\section{Maps preserving triple transition pseudo-probabilities between minimal tripotents}

As we recalled at the introduction, the \emph{transition probability} between two minimal projections $p= \xi\otimes \xi$ and $q= \eta\otimes \eta$ in $B(H)$ is given by $\hbox{tr}(p q) = \hbox{tr}(pq^*) = \hbox{tr}(q p^*) = |\langle \xi, \eta\rangle|^2.$ Let us observe that each minimal projection $p= \xi\otimes \xi$ in $B(H)$ is bi-univocally associated with a pure normal state $\varphi_p\in B(H)_*$ (i.e. an extreme point of the normal state space) at which $p$ attains its norm. Clearly $\varphi_p$ is identified with the pure normal state given by $\varphi (a) = (\xi\otimes \xi) (a) := \langle a(\xi), \xi\rangle = \hbox{tr} ( a p )$ ($a\in B(H)$). Thus, the transition probability between $p$ and $q$ is given by the identity \begin{equation}\label{eq triple trans probability in terms of the evaluation at the supporting pure state} \hbox{tr}(p q) = |\langle \xi, \eta\rangle|^2 = |\varphi_p(q)|^2 = |\varphi_q (p)|^2.
\end{equation}

For each minimal partial isometry $e = \xi \otimes \eta$ in $B(H)$, with $\xi,\eta$ unitary vectors in $H$, there exists a unique extreme point $\varphi_e$ of the closed unit ball of $C_1(H) = B(H)_*$ such that $\varphi_e (e) = 1$. Actually $\varphi_e$ is defined by $\varphi_e (x) := \langle x(\xi), \eta\rangle = \hbox{tr} (e^* x)$ ($x\in B(H)$). Motivated by the identity in \eqref{eq triple trans probability in terms of the evaluation at the supporting pure state}, for each couple $e,v$ of minimal partial isometries in $B(H)$, we define the \emph{triple transition pseudo-probability} between $e$ and $v$ as the scalar $\varphi_e (v)$ --this is not a real probability, since it actually takes complex values. The question is whether we can extend this definition to the wider setting of Cartan factors and atomic JBW$^*$-triples?\smallskip

The lacking of a positive cone in general JB$^*$-triples induced us to replace the lattice of projections in $B(H)$ by the poset of tripotents in a Cartan factor or in an atomic JBW$^*$-triple in our recent study on bijections preserving the partial ordering and orthogonality between the poset of two atomic JBW$^*$-triples in \cite{FriPe2021}. Here we introduce the {triple transition pseudo-probability} between two minimal tripotents in an atomic JBW$^*$-triple. To understand well the definition we need to recall some geometric properties of JBW$^*$-triples. Following \cite{FriRu85}, the extreme points of the closed unit ball, $\mathcal{B}_{M_*},$ of the predual, $M_*$, of a JBW$^*$-triple $M$ are called \emph{atoms} or \emph{pure atoms}. We recall that the extreme points of the convex set of all positive functionals with norm $\leq 1$ in the predual of a von Neumann algebra are called \emph{pure states}. The symbol $\partial_{e} (\mathcal{B}_{M_*})$ will stand for the set of all pure atoms of $M$.\smallskip

By \cite[Proposition 4]{FriRu85} for each minimal tripotent $e$ in a JBW$^*$-triple $M$ there exists a unique pure atom $\varphi_e$ satisfying $P_2 (e) (x) = \varphi_e(x) e$ for all $x\in M$. Furthermore, the mapping $$\mathcal{U}_{min} (M)\to \partial_{e} (\mathcal{B}_{M_*}), \ \  e\mapsto \varphi_e $$ is a bijection from the set of minimal tripotents in $M$ onto the set of pure atoms of $M$. \smallskip

We are now in a position to introduce the key notion of this note. 

\begin{definition}\label{def triple transition pseudo-probability} Let $e$ and $v$ be minimal tripotents in a JBW$^*$-triple $M$. We define the \emph{triple transition pseudo-probability} from $e$ to $v$ as the complex number given by 
	\begin{equation}\label{TTPp}
		TTP(e,v)= \varphi_v(e).
	\end{equation}
\end{definition}

Observe that every triple transition pseudo-probability lies in the closed unit ball of $\mathbb{C}$. Formally speaking, the triple transition pseudo-probability is not a probability because it can take complex values. However, it satisfies many interesting and natural properties. For example, by \cite[Lemma 2.2]{FriRu85} we have \begin{equation}\label{eq pseudo trans prob} TTP(v,e) = \varphi_e (v) = \overline{\varphi_v (e)} = \overline{TTP(e,v)},
\end{equation} for every $e,v\in \mathcal{U}_{min} (M)$, which is naturally expressing the property of symmetry of the triple transition pseudo-probability.\smallskip

If $p$ and $q$ are two minimal projections in a von Neumann algebra $W$, having in mind that $\varphi_p$ is a norm-one functional attaining its norm at $p$, it follows that $\varphi_p$ is a positive normal state on $W$, and hence $TTP (q,p) = \varphi_p (q)$ is a real number in the interval $[0,1]$ and coincides with $TTP (p,q) = \varphi_q (p)$. Therefore the new notion of triple transition pseudo-probability agrees with the usual transition probability in the case of minimal projections. \smallskip   

Moln{\'a}r's theorem \cite[Theorem 2]{Molnar2002}, presented as Theorem \ref{t Molnar minimal pi} in the introduction, can be now restated in the following terms: Let $\Phi: \mathcal{U}_{min} (B(H))\to \mathcal{U}_{min} (B(H))$ be a bijective mapping preserving triple transition pseudo-probabilities. Then $\Phi$ extends to a surjective complex-linear isometry. Inspired by Moln{\'a}r's result, it seems natural to study the bijections preserving the triple transition pseudo-probabilities between the sets of minimal tripotents of two atomic JBW$^*$-triples. The first unexpected conclusion appears when dealing with rank-one JB$^*$-triples. Contrary to the serious obstacles affecting bijective preservers of partial ordering in both directions and orthogonality in the case of rank-one Cartan factors cf. \cite[Remark 3.6]{FriPe2021}), preservers of triple transition pseudo-probabilities between sets of minimal tripotents have an excellent behaviour in the case of rank-one Cartan factors.\smallskip

Let us first recall that a subset $\mathcal{S}$ of a JB$^*$-triple $E$ is called \emph{orthogonal} if $0\notin \mathcal{S}$ and $a\perp b$ for all $a,b\in \mathcal{S}$. The minimal cardinal number $r$ satisfying $\hbox{card}(S) \leq r$ for every orthogonal subset $S \subseteq E$ is called the \emph{rank} of $E$. Spin factors have rank 2 and the exceptional Cartan factors of type 5 and 6 have ranks 2 and 3,  respectively. A JB$^*$-triple has finite rank if and only if it is reflexive (cf. \cite[Proposition 4.5]{BuChu92} and \cite[Theorem 6]{ChuIo90} or \cite{BeLoRo03, BeLoPeRo04}). Furthermore, if $E$ is a JB$^*$-triple of rank-one, it must be reflexive and a rank-one Cartan factor, and moreover, it must be isometrically isomorphic to a complex Hilbert space (see the discussion in \cite[\S 3]{BeLoPeRo04} and  \cite[Table 1 in page 210]{Ka97}).\smallskip

The rank of a tripotent $e$ in a JB$^*$-triple $E$ is defined as the rank of $E_2(e)$. It is known that for each tripotent $e$ in a Cartan factor $C$ we have $r(e) = r(C_2(e)) = n<\infty$ if and only if it can be written as an orthogonal sum of $n$ mutually orthogonal minimal tripotents in $C$ (see, for example, \cite[page 200]{Ka97}).\smallskip

The rank theory plays a fundamental role in the different solutions to Tingley's problem in the case of compact C$^*$-algebras \cite{PeTan19} and weakly compact JB$^*$-triples \cite{FerPe17c,FerPe18Adv}, as well as to prove that every JBW$^*$-triple satisfies the Mazur--Ulam property \cite{BeCuFerPe2021, KalPe2021}. In our next result we shall apply some of the techniques developed in the just quoted results. 

\begin{proposition}\label{p rank-one Cartan factors} Let $\Phi : \mathcal{U}_{min}(E) \to \mathcal{U}_{min}(F)$ be a transformation preserving triple transition pseudo-probabilities, that is, $$TTP(\Phi(u),\Phi(e))= \varphi_{\Phi(e)} (\Phi(u)) = \varphi_{e} (u) = TTP(u,e), \hbox{ for all } e,u\in \mathcal{U}_{min} (E),$$ where $E$ and $F$ are two rank-one JB$^*$-triples. Then $\Phi$ extends to a {\rm(}complex-{\rm)}linear isometric triple homomorphism from $E$ to $F$.
\end{proposition}  

\begin{proof} As we have seen before the statement of this proposition, we can assume that $E$ and $F$ are two complex Hilbert spaces regarded as type 1 Cartan factors. We observe that $\mathcal{U} (E)\backslash \{0\} = \mathcal{U}_{min} (E) = S_{_E}$, the unit sphere of $E,$ and $\mathcal{U} (F)\backslash \{0\} = \mathcal{U}_{min} (F) = S_{_F}$. Since for each $e\in S_{_E},$ $\varphi_e$ is precisely the functional given by $\varphi_e (x) = \langle x, e\rangle$ {\rm(}$x\in E${\rm)}, the hypothesis on $\Phi$ is equivalent to $$\langle  \Phi(u),  \Phi(e)\rangle = \langle u,e\rangle, \hbox{ for all } e,u\in \mathcal{U}_{min} (E) = S_{_E}.$$ A simple computation shows that $$\begin{aligned}
		\|\Phi(e)-\Phi(v)\|^2 &= \langle \Phi(e)-\Phi(v) , \Phi(e)-\Phi(v)\rangle  \\
		&= \langle \Phi(e), \Phi(e)\rangle - \langle \Phi(v) , \Phi(e)\rangle - \langle \Phi(e) , \Phi(v)\rangle + \langle \Phi(v) , \Phi(v)\rangle   \\
		&= \langle e, e\rangle - \langle v , e\rangle - \langle e , v\rangle + \langle v , v\rangle = \| e-v\|^2,
	\end{aligned}$$ for all $e,v\in S(E)$. That is $\Phi : S_{_E} \to S_{_F}$ is an isometry. Moreover, by the assumptions on $\Phi$ we also have $$\langle -\Phi(e), \Phi(-e) \rangle = -\langle \Phi(e), \Phi(-e) \rangle = - \langle e, -e \rangle = 1,$$ which proves that $\Phi(-e) = -\Phi(e)$, for all $e\in S_{_E}$. An application of the solution to Tingley's problem for Hilbert spaces established by G.G. Ding in \cite[Theorem 2.2]{Ding2002} guarantees the existence of a real linear isometry $T: E\to F$ whose restriction to $S_{_E}$ is $\Phi$.\smallskip

We shall finally show that $T$ is complex linear. As before, by the assumptions on $\Phi$, for each $\lambda\in \mathbb{T}$ we also have $$\langle \lambda \Phi(e), \Phi(\lambda e) \rangle = \lambda \langle \Phi(e), \Phi(\lambda e) \rangle = \lambda \langle e, \lambda e \rangle = 1,$$ witnessing that $\Phi(\lambda e) = \lambda \Phi(e)$, for all $e\in S_{_E}$ and $\lambda\in \mathbb{T}$. The rest is clear.   
\end{proof}

Let us note that in the previous proposition we are not assuming that $\Phi$ is injective nor surjective. \smallskip

 It is now time to see a handicap or a limitation of the triple transition pseudo-probability. Let $p$ and $q$ be two minimal projections in a von Neumann algebra $W$. Suppose that the transition probability between $p$ and $q,$ as given in \eqref{eq triple trans probability in terms of the evaluation at the supporting pure state}, is zero, that is $\varphi_p (q) =0$, or equivalently, $P_2 (p) (q) = p q p =0$. Since $0= p q p = (pq)(pq)^*$, we deduce that $p q = q p =0,$ and thus $ q = (1-p) q (1-q) = P_0(p) (q) \perp p$. This property does not always hold when projections are replaced with tripotents or partial isometries, for example if $e$ and $v$ are minimal tripotents in a Cartan factor $C$ with $v\in C_1(e)$ we clearly have $\varphi_e (v) =0$ but $e$ and $v$ are not orthogonal. A simple example can be given by $e = \left(                                                                             \begin{array}{cc} 1 & 0 \\
 	0 & 0 \\
  \end{array}
\right)$ and $v = \left(
\begin{array}{cc}
	 0 & 1 \\
	  0 & 0 \\
  \end{array}
 \right)$ in $M_2(\mathbb{C})$. However, every (real linear) triple homomorphism between JB$^*$-triples preserves orthogonality. Despite of this handicap, we can now get a first extension of every bijective transformation preserving triple transition probabilities between the sets of minimal tripotents of two atomic JBW$^*$-triples to their socles.\smallskip
 
Let us first recall some structure results for atomic JBW$^*$-triples. Every JBW$^*$-triple $M$ decomposes as the orthogonal sum of two weak$^*$-closed ideals $\mathcal{A}$ and $\mathcal{N}$, where $\mathcal{A}$ is an atomic JBW$^*$-triple (called the \emph{atomic part} of $M$) and $\mathcal{N}$ contains no minimal tripotents \cite[Theorem 2]{FriRu85}. Furthermore, $M_*$ decomposes as the $\ell_1$-sum of two norm closed subspaces $\mathcal{A}_{*}$--the predual of $\mathcal{A}$-- and $\mathcal{N}_{*}$--the predual of $\mathcal{N}$-- satisfying that $\mathcal{A}_{*}$ is the norm closure of the linear span of all pure atoms of $M$ and the closed unit ball of $\mathcal{N}_{*}$ contains no extreme points \cite[Theorem 1]{FriRu85}.\smallskip

At this stage the reader should also get some information about elementary JB$^*$-triples. Let $C_j$ be  a Cartan factor of type $j\in \{1,\ldots, 6\}.$ The elementary JB$^*$-triple, $K_j,$ of type $j$ associated with $C_j$ is defined as follows: $K_1 = K (H_1, H_2)$; $K_i = C \cap K(H)$ when $C$ is of type $i = 2 , 3$, and $K_j = C_j$ in the remaining cases (cf. \cite{BuChu92}). For each elementary JB$^*$-triple of type $j$, its bidual space is precisely a Cartan factor of $j$.\smallskip
 
The socle of a JB$^*$-triple $E$, soc$(E)$, is the (non-necessarily closed) linear subspace of $E$ generated by all minimal tripotents in $E$. For example, the socle of $B(H)$ is the subspace, $\mathcal{F} (H)$, of all finite rank operators, and it is not, in general, closed. If $C$ is a Cartan factor of finite rank (or, more generally, a reflexive JB$^*$-triple), every element in $C$ can be written as a finite linear combination of mutually orthogonal minimal tripotents (see \cite[Proposition 4.5 and Remark 4.6]{BuChu92} or \cite{BeLoPeRo04}), and thus the socle of $C$ is the whole $C$--that is soc$(C)= \mathcal{K}(C) = C$. For a general Cartan factor we have $\overline{\hbox{soc}(C)}^{\|-\|} = \mathcal{K} (C)$ and $\overline{\mathcal{K} (C)}^{w^*} = C$. In an atomic JBW$^*$-triple $M$, the symbol $\mathcal{K} (M)$ will stand for the $c_0$-sum of the elementary JB$^*$-triples associated with the Cartan factors expressing $M$ as an $\ell_{\infty}$-sum.

\begin{theorem}\label{t bijections preserving triple transition pseudo-probabilities} Let $\Phi : \mathcal{U}_{min}(M) \to \mathcal{U}_{min}(N)$ be a bijective transformation preserving triple transition pseudo-probabilities {\rm(}i.e., $TTP(\Phi(v), \Phi(e))=\varphi_{\Phi(e)} (\Phi(v)) = \varphi_{e} (v)=TTP(v,e),$ for all $e,v$ in $\mathcal{U}_{min} (M)${\rm)}, where $M$ and $N$ are atomic JBW$^*$-triples. Then there exists a bijective {\rm(}complex{\rm)} linear mapping $T_0: \hbox{soc}(M) \to \hbox{soc}(N)$ whose restriction to $\mathcal{U}_{min} (M)$ is $\Phi$. 
\end{theorem}

\begin{proof}[Proof of Theorem \ref{t bijections preserving triple transition pseudo-probabilities}] Clearly, the pure atoms of $M$ and $N$ are norming sets for $\mathcal{K}(M)$ and $\mathcal{K}(N)$, respectively. Let us suppose that $\displaystyle \sum_{i = 1}^m \alpha_i e_i = \sum_{j=1}^{m} \beta_j v_j \in \hbox{soc}(M),$ where $\alpha_i, \beta_j\in \mathbb{C}$ and $e_i, v_j\in \mathcal{U}_{min} (M)$. By the hypothesis on $\Phi$, for each $\psi \in \partial_{e} (\mathcal{B}_{M_*})$, there exists $\Phi(w)= \tilde{w}\in \mathcal{U}_{min} (N)$ (and $w\in \mathcal{U}_{min} (M)$) such that $\psi = \psi_{\tilde{w}}= \psi_{\Phi(w)}$. It also follows from the hypotheses that $$\begin{aligned}\psi\left(\sum_{i = 1}^m \alpha_i \Phi(e_i) \right) & = \psi_{\Phi(w)} \left(\sum_{i = 1}^m \alpha_i \Phi(e_i) \right) = \sum_{i = 1}^m \alpha_i \psi_{\Phi(w)} \left(\Phi(e_i) \right) \\
		&= \sum_{i = 1}^m \alpha_i \psi_{w} \left( e_i \right) = \psi_{w} \left( \sum_{i = 1}^m \alpha_i  e_i \right) = \psi_{w} \left(  \sum_{j=1}^{m} \beta_j v_j  \right)\\
		&=  \sum_{j=1}^{m} \beta_j  \psi_{w} \left( v_j  \right)  = \sum_{j=1}^{m} \beta_j  \psi_{\Phi(w)} \left( \Phi(v_j)  \right) \\
		&= \psi_{\Phi(w)} \left( \sum_{j=1}^{m} \beta_j   \Phi(v_j)  \right)
		= \psi\left(\sum_{j=1}^{m} \beta_j \Phi(v_j)\right).
	\end{aligned}$$ The arbitrariness of $\psi \in \partial_{e} (\mathcal{B}_{M_*})$ together with the fact that the set of pure atoms of $N$ separates the point of $\mathcal{K} (M)$ imply that $$\sum_{i = 1}^m \alpha_i \Phi(e_i) = \sum_{j=1}^{m} \beta_j   \Phi(v_j).$$ Therefore, the mapping $T_0: \hbox{soc} (M)\to \hbox{soc} (N)$, $T_0\Big(\sum_{i = 1}^m \alpha_i e_i\Big) = \sum_{i = 1}^m \alpha_i \Phi(e_i)$ is well-defined and linear. We further know that $T_0 (e) = \Phi (e)$ for all $e\in \mathcal{U}_{min}(M)$.\smallskip
	
	We can similarly define a linear mapping $R_0: \hbox{soc} (N)\to \hbox{soc} (M)$ satisfying $R_0 (\Phi (e)) = e$ for all $e\in \mathcal{U}_{min}(M)$ and $R_0 = T_0^{-1}$. Therefore $T_0$ and $R_0$ are bijections. 
\end{proof}

It should be remarked that, at this stage the hypotheses of the previous Theorem \ref{t bijections preserving triple transition pseudo-probabilities} do not imply, in a simple way, that the linear mapping $T_0$ is continuous. Actually, if $e_1, \ldots, e_n$ are mutually orthogonal minimal tripotents in $M$, we have $\left\|T_0 \Big(\sum_{j=1}^n e_j \Big)\right\| = \left\| \sum_{j=1}^n \Phi(e_j) \right\|\leq n.$ We cannot get a better bound without assuming orthogonality (and hence $M$-orthogonality) on the minimal tripotents $\Phi(e_1), \ldots, \Phi(e_n).$ In this line we recall next a result by F.J. Herves and J.M. Isidro from \cite{HervIs92}.

\begin{theorem}\label{t HervesIsidro}{\rm\cite[Theorem in page 199]{HervIs92}} Let $E$ be a finite-rank JB$^*$-triple, and let $T: E\to E$ be a	linear mapping (continuity is not assumed). Then the following statements are equivalent:\begin{enumerate}[$(1)$]
\item $T$ is a triple automorphism. 
\item $T\left(\mathcal{U}_{min}(E)\right) = \mathcal{U}_{min}(E)$ and preserves orthogonality.
\end{enumerate}
\end{theorem}

We establish now a hybrid version of the previous two results. 

\begin{corollary}\label{c bijections preserving triple transition pseudo-probabilities and orthogonality} Let $\Phi : \mathcal{U}_{min}(M) \to \mathcal{U}_{min}(N)$ be a bijective transformation preserving orthogonality and triple transition pseudo-probabilities {\rm(}i.e. $TTP(\Phi(v),\Phi(e))=\varphi_{\Phi(e)} (\Phi(v)) = \varphi_{e} (v)=TTP(v,e),$ for all $e,v$ in $\mathcal{U}_{min} (M)${\rm)}, where $M$ and $N$ are atomic JBW$^*$-triples. Then $\Phi$ extends {\rm(}uniquely{\rm)} to a surjective complex-linear {\rm(}isometric{\rm)} triple isomorphism from $M$ onto $N$.
\end{corollary}

\begin{proof} By Theorem \ref{t bijections preserving triple transition pseudo-probabilities} there exists a linear bijection $T_0: soc(M)\to soc(N)$ whose restriction to $\mathcal{U}_{min}(M)$ coincides with $\Phi$. By hypotheses, given $u,v\in \mathcal{U}_{min} (M)$ with $u\perp v$, we have $\Phi (u) \perp \Phi(v)$. Having in mind that for each $x$ in the closed unit ball of $\hbox{soc}(M)$ there exists a finite family $\{e_n\}_n$ of mutually orthogonal minimal tripotents in $M$ and $\{\lambda_n\}_n$ in $\mathbb{R}^+$ such that $\displaystyle x = \sum_n \lambda_n e_n$ and $1 = \|x\| = \max\{\lambda_n : n\}$ (cf. \cite[Remark 4.6]{BuChu92}). It follows from the definition of $T_0$ that $$T_0(x) = T_0 \left( \sum_n \lambda_n e_n \right) = \sum_n \lambda_n T_0(e_n) = \sum_n \lambda_n \Phi(e_n),$$ and hence $\|T_0(x)\|= \|x\|,$ because, by hypotheses, $\{\Phi(e_n)\}_n$ is a family of mutually orthogonal minimal tripotents in $N$. Furthermore, by the previous conclusion $$\begin{aligned} \{T_0(x), T_0(x), T_0(x)\} & = \sum_n \lambda_n^3 \{\Phi(e_n),\Phi(e_n),\Phi(e_n)\} \\
&= \sum_n \lambda_n^3 \Phi(e_n) = \sum_n \lambda_n^3 T_0(e_n) = T_0(\{x,x,x\}),
\end{aligned}$$ which shows that $T_0$ is a contractive triple isomorphism from soc$(M)$ onto soc$(N)$.\smallskip

We can therefore find a continuous linear extension of $T_0$ to a continuous (isometric) linear triple isomorphism from $\mathcal{K} (M)$ onto $\mathcal{K} (N)$ denoted by the same symbol $T_0$. The bitransposed mapping $T_0^{**}:  \mathcal{K} (M)^{**} = M \to \mathcal{K} (N)^{**} =N$ is a triple isomorphism whose restriction to $\mathcal{U}_{min} (M)$ coincides with $\Phi$. This finishes the proof of the result.
\end{proof}

It seems a natural (and important) question to ask whether a bijection preserving triple transition pseudo-probabilities between the sets of minimal tripotents in two atomic JBW$^*$-triples also preserves orthogonality. That is, whether in Corollary \ref{c bijections preserving triple transition pseudo-probabilities and orthogonality} the hypothesis concerning preservation of orthogonality can be relaxed. This will be answered for spin and type 1 Cartan factors along the next sections. \smallskip

\section{The case of spin factors}\label{sec: spin factors}

As well as the study of those maps preserving triple transition pseudo-probabilities between the sets of minimal tripotents in two rank-one JB$^*$-triples deserved its own treatment in Proposition \ref{p rank-one Cartan factors}, the case of spin factors is also worth to study by itself. \smallskip

Let us fix a spin factor $M$ whose inner product, involution and triple product are given by $\langle \cdot, \cdot \rangle,$ $x\mapsto \overline{x}$, and  $$\{a,b,c\} = \langle a, b\rangle c + \langle c,b\rangle a-\langle a,\bar{ c}\rangle \bar{b},$$ respectively (cf. the definition in page \pageref{def spin factor}). It is usually assumed that dim$(M)\geq 3$; actually if dim$(M)=2$, the defined structure produces $\mathbb{C}\oplus^{\infty}\mathbb{C}$, which is not a factor (cf. \cite[Remark 4.3]{Ka97}). The real subspace $$M_{_\mathbb{R}}^{-}=\{a\in M: \, a=\bar{a}\},$$ of all fixed points for the involution $\overline{\cdot}$ is a real Hilbert space with respect to the restricted inner product $\langle a, b\rangle = \Re\hbox{e}\langle a, b\rangle$ ($a,b\in M_{_\mathbb{R}}^{-}$), and $M = M_{_\mathbb{R}}^{-}\oplus i M_{_\mathbb{R}}^{-}$. We shall also make use of the real Hilbert space $\mathcal{H}= M_{_\mathbb{R}}$ given by the real underlying space of $M$ equipped with the inner product $\Re\hbox{e} \langle .,.\rangle$. Clearly $M_{_\mathbb{R}}^{-}$ is a closed subspace of $M_{_\mathbb{R}}$. The symbol $\perp_2$ will denote orthogonality in the Euclidean sense.  \smallskip

Each triple automorphism $\Phi$ on the spin factor $M$ is precisely described in the following form:  $$T(a +ib) =\lambda (U(a) + i U(b)) \hbox{ for all $a, b\in M_{_\mathbb{R}}^{-}$,}$$ where $\lambda\in \mathbb{T}$ and $U: M_{_\mathbb{R}}^{-}\to M_{_\mathbb{R}}^{-}$ is a unitary operator (cf. \cite[Theorem in page 196]{HervIs92} or \cite[Section 3.1.3]{Fri2005}).\smallskip

The set of tripotents in $M$ has been intensively studied along the last forty years. If we exclude the zero tripotent, $M$ only contains tripotents of rank-one (minimal) and of rank-two (maximal and unitaries). In the second case we have $$\mathcal{U}_{max} (M) = \{ \lambda a : \lambda\in \mathbb{T}, \ a\in S_{_{M_{_\mathbb{R}}^{-}}} \},$$ while \begin{equation}\label{eq min trip in spin} \mathcal{U}_{min} (M) = \left\{ \frac{a+i b}{2} : a, b\in S_{_{M_{_\mathbb{R}}^{-}}}, \ \langle a, b\rangle =0 \right\}
\end{equation} (see, for example, \cite[Section 3.1.4]{Fri2005} or \cite[Lemma 6.1]{HamKalPe20}, \cite{FriPe2021} or \cite[Section 3]{KalPe2021}).\smallskip

It is well-known, and easy to check, that for each minimal tripotent $v = \frac{a+i b}{2}$ with $a, b\in S_{_{M_{_\mathbb{R}}^{-}}},$  $\langle a, b\rangle =0$ the Peirce-$0$ and $1$ subspaces are the following:
\begin{equation}\label{eq Peirce 0 and 1 subspace spin} M_0 (v) = \mathbb{C} \overline{v}, \hbox{ and } M_1 (v) = \{x\in M : x\perp_2 v, \overline{v} \}= \{x\in M : x\perp_2 a, b \}.
\end{equation} 

Another important fact to have in mind is the following: for each minimal tripotent $v$ in a spin factor $M$ we have \begin{equation}\label{pure states in spin factors} P_2(v) (x) = 2 \langle x, v\rangle v, \hbox{ and hence } \varphi_v (x) = 2 \langle x, v\rangle\ \ \ (x\in M).
\end{equation}

The next technical lemma is presented separately to simplify the arguments. 

\begin{lemma}\label{l min trip sharing Peirce-1 spin} Let $v$ and $w$ be minimal tripotents in a spin factor $M$. Suppose that $M_1(v) = M_1(w)$. Then $w$ lies in he linear span of $v$ and $\overline{v}$. If we additionally assume that $TTP(w,v) = \varphi_v (w)= 0$, then $w$ belongs to $\mathbb{T} \overline{v} \subset M_0(v).$
\end{lemma}

\begin{proof} Let us write $v = \frac{a+i b}{2}$, $w = \frac{c+i d}{2}$ with $a, b,c,d\in S_{_{M_{_\mathbb{R}}^{-}}}$ and  $\langle a, b\rangle =0 = \langle c, d\rangle$. Since $\{x\in M : x\perp_2 a, b \} = M_1 (v) = M_1 (w) =\{x\in M : x\perp_2 c, d \},$ it follows tha the Euclidean orthogonal complements of the sets $\{a,b\}$ and $\{c,d\}$ in $M$ coincide. It is straightforward to check that, under these circumstances, $c = \alpha_1 a + \alpha_2 b$ and  $d = \beta_1 a + \beta_2 b$, with $(\alpha_1, \alpha_2), (\beta_1,\beta_2)\in S_{\ell_2^2(\mathbb{R})}\equiv \mathbb{T}$ . Therefore $$w = \frac{c+i d}{2} = \frac{\alpha_1 a + \alpha_2 b +i(\beta_1 a + \beta_2 b)}{2} \in \hbox{span} \{v,\overline{v}\},$$ as desired.\smallskip
	
Finally, if $TTP(w,v) = \varphi_v (w)= 0$, the expression of the pure states given in \eqref{pure states in spin factors} assures that $$ 0 = \varphi_v(w) = 2 \langle w, v\rangle = 2 \left\langle \frac{\alpha_1 a + \alpha_2 b +i(\beta_1 a + \beta_2 b)}{2}, \frac{a+i b}{2}\right\rangle,$$ or equivalently, 
$$ \alpha_1 +\beta_2 =0 = \beta_1 -\alpha_2,$$ that is, $\beta_1 + i \beta_2 = \alpha_2 - i \alpha_1 = (-i) (\alpha_1 + i \alpha_2),$ and thus $$w =  \frac{(\alpha_1 +i \alpha_2) a + (\alpha_2 -i \alpha_1) b}{2} = (\alpha_1 +i \alpha_2)\frac{a-i b}{2} \in \mathbb{T} \overline{v}, $$ which finishes the proof.
\end{proof}

We can now prove that every bijection preserving triple transition pseudo-proba-bilities between the sets of minimal tripotents in two spin factors also preserves orthogonality, and consequently extends to a triple isomorphism between the spin factors.  

\begin{theorem}\label{t bijections preserving triple transition pseudo-probabilities spin factors} Let $\Phi : \mathcal{U}_{min}(M) \to \mathcal{U}_{min}(N)$ be a bijective transformation preserving triple transition pseudo-probabilities {\rm(}i.e., $TTP(\Phi(v), \Phi(e))=\varphi_{\Phi(e)} (\Phi(v)) = \varphi_{e} (v) = TTP(v,e),$ for all $e,v$ in $\mathcal{U}_{min} (M)${\rm)}, where $M$ and $N$ are spin factors. Then there exists a {\rm(}unique{\rm)} triple isomorphism $T: M \to N$ whose restriction to $\mathcal{U}_{min} (M)$ is $\Phi$. 
\end{theorem}

\begin{proof} We shall denote by the same symbols $\langle.,.\rangle$ and $x\mapsto \overline{x}$ the inner products and the involutions on $M$ and $N$. Let $T_0: \hbox{soc}(M) \to \hbox{soc}(N)$ be the linear bijection  whose restriction to $\mathcal{U}_{min} (M)$ is $\Phi$ given by Theorem \ref{t bijections preserving triple transition pseudo-probabilities}. Since $M$ and $N$ are spin factors with rank-two we have soc$(M) = M$ and soc$(N) = N.$ Furthermore, every element $x$ in $M$ writes in the form $x = \lambda_1 v_1 + \lambda_2 v_2$, where $v_1$ and $v_2$ are two orthogonal minimal tripotents in $M$, $\lambda_i\in [0,\|x\|]$. In particular, $\|T_0(x)\| \leq \lambda_1 \|\Phi(v_1)\| + \lambda_2 \|\Phi(v_2)\|\leq 2 \|x\|$. Therefore $T_0$ is a bounded linear bijection from $M$ onto $N$ sending minimal tripotents to minimal tripotents and preserving triple transition pseudo-probabilities.\smallskip

Let us take $a,b\in S_{_{M_{_\mathbb{R}}^{-}}}$ with $\langle a, b \rangle =0$. Let us write $T_0(a) = a_1 + i a_2$ and $T_0(b) = b_1 + i b_2$, with $a_j,b_j\in  N_{_\mathbb{R}}^{-}.$ Since the elements $T_0( \frac{a+i b}{2}) =  \frac{(a_1-b_2)+i (b_1+a_2)}{2}$ and $T_0( \frac{a-i b}{2}) =  \frac{(a_1+b_2)+i (a_2-b_1)}{2}$ are minimal tripotents we deduce that the following identities hold: $$ \left\{\begin{array}{l} \| a_1- b_2\|_2^2 =\| a_2+b_1\|_2^2 =1,\\
	\| a_1+b_2\|_2^2 =\| a_2-b_1\|_2^2 =1,\\
\langle a_1+b_2,  a_2-b_1\rangle =0,\\
\langle a_1- b_2 ,  a_2+b_1\rangle =0, 
\end{array} \right. ,$$ equivalently, \begin{equation}\label{eq second list of identities from min trip}  \left\{\begin{array}{l} \| a_1\|_2^2 + \| b_2\|_2^2 =\| b_1\|_2^2 +\| a_2\|_2^2 =1,\\
	\langle a_1, b_2\rangle =0 = \langle a_2, b_1\rangle, \\
	\langle a_1, a_2\rangle - \langle b_1, b_2\rangle =0, \\
	\langle a_1, b_1\rangle - \langle a_2, b_2\rangle =0.
\end{array} \right.
\end{equation} 	

On the other hand, by \eqref{pure states in spin factors}, for each minimal tripotent $v$ in a spin factor we have $\varphi_v (x) = 2 \langle x, v\rangle$. Applying now that $\Phi$ (and $T_0$) preserves triple transition pseudo-probabilities among minimal tripotents and the explicit expression of the pure atoms given above we have
$$\begin{aligned} 0= 2 \left\langle\frac{a+i b}{2}  ,  \frac{a-i b}{2} \right\rangle &=\varphi_{\frac{a-i b}{2}} \left( \frac{a+i b}{2} \right) = \varphi_{T_0\left( \frac{a-i b}{2}\right)} \left( T_0\left( \frac{a+i b}{2}\right) \right) \\
&= 2 \left\langle T_0\left( \frac{a+i b}{2}\right) , T_0\left( \frac{a-i b}{2}\right)\right\rangle \\
&= 2 \left\langle \frac{(a_1-b_2)+i (b_1+a_2)}{2},   \frac{(a_1+b_2)+i (a_2-b_1)}{2}\right\rangle	
\end{aligned}.$$ By computing the imaginary parts in the previous equality we arrive at 
$$\begin{aligned} 0 &= \left\langle \frac{ (b_1+a_2)}{2},   \frac{(a_1+b_2)}{2}\right\rangle - \left\langle \frac{(a_1-b_2)}{2},   \frac{ (a_2-b_1)}{2}\right\rangle	
\end{aligned}$$ equivalently, $$\langle a_1, b_1\rangle + \langle a_2, b_2\rangle =0,$$ which combined with the equality in the fourth line of \eqref{eq second list of identities from min trip} gives $$\langle a_1, b_1\rangle =  \langle a_2, b_2\rangle =0.$$ The last identity together with the conclusion in the second line of  \eqref{eq second list of identities from min trip} show that $$\begin{aligned}&\langle a_1, b_1\rangle =  \langle a_2, b_2\rangle = \langle a_1, b_2\rangle = \langle a_2, b_1\rangle=0 \\ 
&\Leftrightarrow a_j\perp_2 b_k, \ j,k = 1,2\hbox{ in the Hilbert space }\! N_{\mathbb{R}}.
\end{aligned}$$ We have therefore shown that \begin{equation}\label{eq T0 preserves Hilbert orthogonality} \boxed{a\!\perp_2\! b \hbox{ in the Hilbert space }\! M_{_\mathbb{R}}^{-}}\! \Rightarrow\! \boxed{T_0(a) \perp_2 T_0 (b) \hbox{ in the Hilbert space } N_{_\mathbb{R}}.\!}  
\end{equation} 

Therefore $T_0|_{M_{_\mathbb{R}}^{-}} : M_{_\mathbb{R}}^{-} \to N_{_\mathbb{R}}$ is an injective linear mapping preserving orthogonality in the Hilbert sense, that is, $\langle a,b\rangle =0\Rightarrow \langle T_0(a), T_0(b)\rangle =0 $. Another interesting result on preservers, proved by J. Chmieli{\'n}ski, assures that  $T_0|_{M_{_\mathbb{R}}^{-}}$ is a positive multiple of an isometry, that is, there exists a positive $\gamma$ and a real linear isometry $U :{M_{_\mathbb{R}}^{-}} \to N_{_\mathbb{R}}$ such that $T_0|_{M_{_\mathbb{R}}^{-}} = \gamma  U$ (cf. \cite[Theorem 1]{Chmielinski05}). \smallskip

We shall next prove that $\gamma=1$. Indeed, let us fix two orthogonal elements $a,b$ in the unit sphere of  ${M_{_\mathbb{R}}^{-}}$. By hypothesis, $\widehat{v} = T_0(\frac{a+i b}{2}) = \frac{\gamma U(a) + i \gamma U(b)}{2} = \gamma \frac{U(a) + i U(b)}{2}$ is a minimal tripotent in $N$. It is well-known that in this case
\begin{equation}\label{eq hatv and overline hatv are orthogonal} \left\langle \frac{U(a) + i  U(b)}{2}, \frac{\overline{ U(a)} - i \overline{ U(b)}}{2} \right\rangle = \left\langle \gamma^{-1} \widehat{v}, \gamma^{-1} \overline{\widehat{v}} \right\rangle = \gamma^{-2} \left\langle \widehat{v}, \overline{\widehat{v}} \right\rangle = 0.
\end{equation} By applying the properties of $U$ and the previous identity we get 
$$\begin{aligned}&\left\{\frac{U(a) + i U(b)}{2}, \frac{U(a) + i U(b)}{2}, \frac{U(a) + i U(b)}{2} \right\}\\
& = 2 \left\langle \frac{U(a) + i  U(b)}{2}, \frac{{ U(a)} + i { U(b)}}{2} \right\rangle \frac{U(a) + i U(b)}{2} \\
&- \left\langle \frac{U(a) + i  U(b)}{2}, \frac{\overline{ U(a)} - i \overline{ U(b)}}{2} \right\rangle  \frac{\overline{ U(a)} - i \overline{ U(b)}}{2} \\
& = \frac{1}{2} \left( \|U(a)\|_2^2 + \|U(b)\|_2^2\right) \frac{U(a) + i U(b)}{2} \\ 
&= \frac{1}{2} \left( \|a\|_2^2 + \|b\|_2^2\right) \frac{U(a) + i U(b)}{2} =  \frac{U(a) + i U(b)}{2}\\
\end{aligned},$$ witnessing that $\frac{U(a) + i U(b)}{2}$ is a tripotent. Since $\frac{U(a) + i U(b)}{2} = \gamma^{-1} \widehat{v}$ we obtain $\gamma =1$ as desired.\smallskip

Let us go back to the identity in \eqref{eq hatv and overline hatv are orthogonal} to deduce that \begin{equation}\label{eq 14}\begin{aligned} 0 &= \left\langle U(a) + i  U(b), \overline{ U(a)} - i \overline{ U(b)} \right\rangle \\ 
		&= \left\langle U(a) , \overline{ U(a)}  \right\rangle - \left\langle U(b),  \overline{ U(b)} \right\rangle + i \left\langle U(a),  \overline{ U(b)} \right\rangle + i \left\langle  U(b), \overline{ U(a)} \right\rangle \\
		&= \left\langle U(a) , \overline{ U(a)}  \right\rangle - \left\langle U(b),  \overline{ U(b)} \right\rangle + 2 i \left\langle U(a),  \overline{ U(b)} \right\rangle. 
\end{aligned}
\end{equation} Since in this argument the roles of $a$ and $b$ are completely symmetric, by replacing $\frac{a +i b}{2}$ with $\frac{b +i a }{2}$ we derive  \begin{equation}\label{eq 15}\begin{aligned} 0 &= \left\langle U(b) , \overline{ U(b)}  \right\rangle - \left\langle U(a),  \overline{ U(a)} \right\rangle + 2 i \left\langle U(b),  \overline{ U(a)} \right\rangle\\
	&= \left\langle U(b) , \overline{ U(b)}  \right\rangle - \left\langle U(a),  \overline{ U(a)} \right\rangle + 2 i \left\langle U(a),  \overline{ U(b)} \right\rangle. 
\end{aligned}
\end{equation} Now, adding \eqref{eq 14} and \eqref{eq 15} we get \begin{equation}\label{eq Ua and overline Ub are orthogonal} \left\langle T_0(a), \overline{T_0(b)}\right\rangle = \left\langle U(a), \overline{U(b)}\right\rangle = 0, \hbox{ for all } a,b\in S_{_{M_{\mathbb{R}}^-}} \hbox{ with } a\perp_2 b.
\end{equation}

We shall next show that for any minimal tripotent $v= \frac{a+i b} {2}$ the identities \begin{equation}\label{eq T0 preserves Peirce-1}  N_1(T_0(\overline{v}))= T_0\left( M_1(\overline{v}) \right) = T_0\left( M_1(v) \right) = N_1(T_0(v)) = N_1(\overline{T_0(v)}) 
\end{equation} hold. The second and fourth equalities are clear because $M_1(\overline{v}) = M_1(v)$ and $N_1(T_0(v))$ $= N_1(\overline{T_0(v)})$ (cf. \eqref{eq Peirce 0 and 1 subspace spin}). For the third one suppose that dim$(M)= \hbox{dim}(M_{\mathbb{R}}^-)\geq 4$, and observe that in this case every minimal tripotent in $M_1(v)$ is of the form $\frac{c + id }{2}$ with $c,d\in M_{\mathbb{R}}^{-}$ and $c,d\perp_2 a,b$, and every element in $M_1(v)$ writes as the linear combination of a minimal tripotent of this form and its orthogonal. Having in mind \eqref{eq Ua and overline Ub are orthogonal} and the properties of $U$ we obtain $$\begin{aligned}\left\{T_0(v), T_0(v), T_0\left(\frac{c+ i d}{2}\right)  \right\} &= \left\{\frac{U(a) + i U(b)}{2}, \frac{U(a) + i U(b)}{2}, \frac{U(c) + i U(d)}{2}\right\} \\
&=\left\langle \frac{U(a) + i U(b)}{2}, \frac{U(a) + i U(b)}{2} \right\rangle \frac{U(c) + i U(d)}{2} \\
&+  \left\langle \frac{U(c) + i U(d)}{2}, \frac{U(a) + i U(b)}{2} \right\rangle \frac{U(a) + i U(b)}{2} \\
&- \left\langle \frac{U(a) + i U(b)}{2}, \frac{\overline{U(c)} - i \overline{U(d)}}{2} \right\rangle \frac{\overline{U(a)} - i \overline{U(b)}}{2} \\
&= \frac{1}{2} \frac{U(c) + i U(d)}{2},
\end{aligned}$$ witnessing that $T_0\left(\frac{c+ i d}{2}\right)=\frac{U(c) + i U(d)}{2}\in N_1\left(T_0(v)\right),$ and consequently $T_0\left( M_1(v) \right)$ $\subseteq N_1(T_0(v))$, but the equality holds by the bijectivity of $T_0$. \smallskip

If dim$(M)= \hbox{dim}(M_{\mathbb{R}}^-) =3,$ then $M_1(v)$ is one dimensional and of the form $\mathbb{C} c$ with $c$ in the unit sphere of $M_{\mathbb{R}}^-$. In this case
$$\begin{aligned}\left\{T_0(v), T_0(v), T_0\left(\alpha c \right)  \right\} &= \left\{\frac{U(a) + i U(b)}{2}, \frac{U(a) + i U(b)}{2}, \alpha U(c)\right\} \\
	&=\left\langle \frac{U(a) + i U(b)}{2}, \frac{U(a) + i U(b)}{2} \right\rangle \alpha U(c) \\
	&+  \left\langle \alpha U(c), \frac{U(a) + i U(b)}{2} \right\rangle \frac{U(a) + i U(b)}{2} \\
	&- \left\langle \frac{U(a) + i U(b)}{2}, \overline{\alpha U(c) } \right\rangle \frac{\overline{U(a)} - i \overline{U(b)}}{2} \\
	&= \frac{1}{2} \frac{U(c) + i U(d)}{2},
\end{aligned}$$ by the properties of $U$ and \eqref{eq Ua and overline Ub are orthogonal}. This shows that $T_0\left( M_1(v) \right)$ $\subseteq N_1(T_0(v))$, and consequently, $T_0\left( M_1(v) \right)$ $= N_1(T_0(v)).$ The same argument shows the validity of the first equality in \eqref{eq T0 preserves Peirce-1}.\smallskip

We have therefore proved that $N_1(T_0(\overline{v}))=  N_1(T_0(v)),$ and since, by hypothesis, $TTP(T_0(\overline{v}), T_0({v})) = TTP(\overline{v},v) =0$, Lemma \ref{l min trip sharing Peirce-1 spin} assures that $T_0(\overline{v})\in \mathbb{T}\ \overline{T_0(v)},$ and thus $T_0(\overline{v}) \perp T_0(v)$ in the spin factor $N$. Since $v = \frac{a+i b}{2}$ is any arbitrary minimal tripotent in $M$ and each minimal tripotent orthogonal to $v$ lies in $\mathbb{T} \overline{v}$, we deduce that $T_0$ preserves orthogonality among $\mathcal{U}_{min} (M)$. Finally, Corollary \ref{c bijections preserving triple transition pseudo-probabilities and orthogonality} asserts that $T_0$ is a triple isomorphism. 
\end{proof}

\section{The case of type 1 Cartan factors}\label{sec: type 1 Cartan factors}

This section is aimed to study those bijections preserving triple transition pseudo-probabilities between the sets of minimal tripotents of two type 1 Cartan factors. In this case we shall try to extend the arguments settled by L. Moln{\'a}r in \cite[Theorem 2]{Molnar2002}. For this purpose we shall focus on those linear operators between type 1 Cartan factors preserving rank-one operators. Let us begin by recalling a result by M. Marcus, B.N. Moyls \cite{MarcusMoyls59} and R. Westwick \cite{Westwick67}.

\begin{theorem}{\rm \cite{MarcusMoyls59, Westwick67}} Let $T: M_{m,n} (\mathbb{C})\to M_{m,n} (\mathbb{C})$ be a linear mapping sending rank-one operators to rank-one operators. Then there exist invertible matrices $u\in M_{m} (\mathbb{C})$ and $M_{n} (\mathbb{C})$ such that one of the next statements holds:
\begin{enumerate}[$(1)$]\item $T (a) = u a v$ for all $a\in M_{m,n} (\mathbb{C})$;
\item $m = n$ and  $T (a) = u a^{t} v$ for all $a\in M_{m,n} (\mathbb{C})$; where $a^t$ denotes the transpose of $a$.
\end{enumerate}
\end{theorem} 

One of the key tools employed by Moln{\'a}r in \cite[Theorem 2]{Molnar2002} is the following consequence of a result due to  M. Omladi\v{c} and P. \v{S}emrl.

\begin{theorem}{\rm \cite[Theorem 3.3]{OmladicSemrl93}} Let $H$ be a complex Hilbert space. Suppose $\Phi : soc(B(H))\to soc(B(H))$ is a surjective linear transformation preserving the rank-one operators in both directions. Then: \begin{enumerate}[$(a)$]\item either there are bijective linear mappings $u,v$ on $H$ such that $\Phi (\xi \otimes \eta)= u(\xi) \otimes v(\eta)$ {\rm(}$\xi,\eta \in H${\rm)};
\item or there are bijective conjugate-linear mappings $u, v$ on $H$ such that
$\Phi (\xi \otimes \eta)=  u(\eta) \otimes v(\xi)$ {\rm(}$\xi,\eta \in H${\rm)}.
\end{enumerate}	
\end{theorem}

It is known that (linear) triple automorphisms on a type 1 Cartan factor of the form $B(H,K)$ are either of the form $T(x) = u x v$ ($x\in B(H,K)$), with $u\in B(K)$ and $v\in B(H)$ unitaries, or dim$(H)=$dim$(K)$ and $T(x) = u x^* v$ ($x\in B(H,K)\equiv B(H)$), with $u,v: H\to H$ anti-unitaries (cf. \cite[page 199]{Ka97}). We establish next the tool required for our purposes. It should be noted that we have strengthened the hypotheses with respect to the mentioned result by Omladi\v{c} and \v{S}emrl, but the current statement will serve for our goals in this note.  

\begin{theorem}\label{t extension of OmladicSemrl} Let $\Phi : soc(B(H_1,K_1))\to soc(B(H_2,K_2))$ be a bijective linear transformation, where $H_1,H_2, K_1$ and $K_2$ are complex Hilbert spaces with dimensions $\geq 2$. Suppose that $\Phi$ preserves rank-one operators in both directions. Then: \begin{enumerate}[$(a)$]\item either there are bijective linear mappings $u:K_1\to K_2,$ and $v:H_1\to H_2$ such that $\Phi (\xi \otimes \eta) =  u(\xi) \otimes v(\eta)$ {\rm(}$\xi\in K_1,\eta \in H_1${\rm)};
\item or there are bijective conjugate-linear mappings $u:H_1\to K_2,$ $v:K_1\to H_2$ such that $\Phi (\xi \otimes \eta)$ $=  u(\eta) \otimes v(\xi)$ {\rm(}$\xi\in K_1,$ $\eta \in H_1${\rm)}.
	\end{enumerate}	
\end{theorem} 

\begin{proof} We shall mimic the notation and arguments by Omladi\v{c} and \v{S}emrl in \cite{OmladicSemrl93}. So, for each $\xi\in K_j,\eta \in H_j$ we set $L_{\xi} :=\{ \xi\otimes \eta : \eta \in H_j\}$ and $R_{\eta}  :=\{ \xi\otimes \eta : \xi \in K_j\}$. It is not hard to check, as in the proof of \cite[Lemma 2.1]{OmladicSemrl93}, that $L_{\xi}$ and  $R_{\eta}$ are maximal among additive subgroups of rank-one operators, that is, every additive group of $soc(H_j,K_j)$ consisting of operators of rank-one is either a subgroup of an $L_{\xi}$, for a vector $\xi\in K$, or a subgroup of an $R_{\eta}$ for $\eta \in H_j$. \smallskip

\emph{Step 1.} We claim that for every $\xi \in K_1$, either there is a $\widehat{\xi}\in K_2$, depending on $\xi$, such that $\Phi( L_{\xi}) = L_{\widehat{\xi}}$, or there is an $\widehat{\eta}\in H_2$, depending on $\xi$, such that $\Phi( L_{\xi}) = R_{\widehat{\eta}}.$ This is clear because $\Phi$ is linear, bijective, and preserves rank-one operators in both directions, and thus it must preserve maximal additive subgroups of rank-one operators. \smallskip

\emph{Step 2.}  If there exists $\xi_0 \in K_1\backslash\{0\}$ such that $\Phi( L_{\xi_0}) = L_{\widehat{\xi_0}}$ for some $\widehat{\xi_0}\in K_2$ depending on $\xi_0$ (respectively, $\Phi( L_{\xi_0}) = R_{\widehat{\eta_0}}$ for some $\widehat{\eta_0}\in H_2$, depending on $\xi_0$), then for each $\xi\in K_1$, $\Phi( L_{\xi}) = L_{\widehat{\xi}}$ for some $\widehat{\xi}\in K_2$ depending on $\xi$ (respectively, $\Phi( L_{\xi}) = R_{\widehat{\eta}}$ for some $\widehat{\eta}\in H_2$, depending on $\xi$). We shall only prove the first statement. Suppose we can find $\xi_0$ and $\xi_1$ in $K_1\backslash\{0\}$ such that $\Phi( L_{\xi_0}) = L_{\widehat{\xi_0}}$ for some $\widehat{\xi_0}\in K_2$ depending on $\xi_0$ and $\Phi( L_{\xi_1}) = R_{\widehat{\eta_1}}$ for some $\widehat{\eta_1}\in H_2$, depending on $\xi_1$. Observing that $L_{\xi_0} = L_{\xi_1}$ if and only if $\xi_0$ and $\xi_1$ are linearly dependent, it follows from the assumptions that  $\xi_0$ and $\xi_1$ must be linearly independent. The element $\widehat{\xi_0}\otimes \widehat{\eta_1}\in \cap L_{\widehat{\xi_0}}\cap  R_{\widehat{\eta_1}} = \Phi( L_{\xi_0}) \cap \Phi( L_{\xi_1})$, and thus we can find $\eta_0,\eta_1\in H_1\backslash\{0\}$ such that $\Phi (\xi_0\otimes \eta_0) = \widehat{\xi_0}\otimes \widehat{\eta_1} = \Phi (\xi_1\otimes \eta_1)$. The injectivity of $\Phi$ assures that $\xi_0\otimes \eta_0 = \xi_1\otimes \eta_1$, however this equality implies that $\xi_0$ and $\xi_1$ are linearly dependent, which is impossible. \smallskip

\emph{Step 3.} Let us assume, by Steps 1 and 2, that for each $\xi \in K_1$, there exists a $\widehat{\xi}\in K_2$, depending on $\xi$, such that $\Phi( L_{\xi}) = L_{\widehat{\xi}}$. Consequently, for each $\eta\in H_1$ there exists a unique $v_{\xi} (\eta)\in H_2$ such that \begin{equation}\label{eq definition property of vxi}  \Phi (\xi\otimes \eta) = \widehat{\xi} \otimes v_{\xi}  (\eta).
\end{equation} It is easy to see that $v_{\xi}$ inherits the linearity of $\Phi$, and so $ v_{\xi}: H_1\to H_2$ is a well-defined linear bijection.\smallskip

Fix a non-zero $\xi_0$ in $K_1$. We shall next prove the existence of a non-zero constant $\tau=\tau(\xi)$ such that $ v_{\xi} = \tau v_{\xi_0}$ for all $\xi \in K_1$.  Namely, fix $\eta\in H_1$. If $\widehat{\xi}$ and $\widehat{\xi_0}$ are linearly independent, let us find, by the hypothesis on $v_{\xi_0},$ and the fact that dim$(H_1)\geq 2$, another element $\eta_1$ such that $v_{\xi_0} (\eta)$ and $v_{\xi_0}(\eta_1)$ are linearly independent. Since the elements $\xi_0 \otimes \eta + \xi \otimes \eta $, $\xi_0 \otimes \eta_1 + \xi \otimes \eta_1 $ and $\xi_0 \otimes (\eta + \eta_1)+ \xi \otimes (\eta + \eta_1)$ have rank-one, the same holds for their images under $\Phi$, that is, the elements $$\left\{\begin{aligned} &\widehat{\xi_0} \otimes v_{\xi_0}(\eta) + \widehat{\xi} \otimes v_{\xi}(\eta),\\ 
&\widehat{\xi_0} \otimes v_{\xi_0}(\eta_1) + \widehat{\xi} \otimes v_{\xi}(\eta_1),\\ 
&\widehat{\xi_0} \otimes v_{\xi_0}(\eta + \eta_1)+ \widehat{\xi} \otimes v_{\xi}(\eta + \eta_1),
\end{aligned}\right.$$ must be rank-one elements. It necessarily follows that each one of the sets $$\{v_{\xi_0}(\eta) ,  v_{\xi}(\eta)\}, \  \{v_{\xi_0}(\eta_1) , v_{\xi}(\eta_1)\}, \ \{v_{\xi_0}(\eta + \eta_1),  v_{\xi}(\eta + \eta_1)\}$$ must be linearly dependent, which implies the existence of a non-zero constant $\alpha$ such that $v_{\xi}(\eta) = \alpha v_{\xi_0}(\eta)$ and $ v_{\xi}(\eta_1) = \alpha v_{\xi_0}(\eta_1)$ for every $\eta$ and $\eta_1$ linearly independent. It follows from the hypothesis on the dimension of $H_1$ that $v_{\xi} = \alpha v_{\xi_0}$. The constant $\alpha$ might depend on the element $\xi$, so we shall denote it by $\alpha(\xi)\in \mathbb{C}$.\smallskip

Suppose now that $\widehat{\xi}$ and $\widehat{\xi_0}$ are linearly dependent. Find another ${\xi_1}$ such that $\widehat{\xi_1}$ and any of $\{\widehat{\xi_0}, \widehat{\xi}\}$ are linearly independent. It follows from the above arguments that $v_{\xi_1} = \alpha (\xi_1) v_{\xi_0},$ $v_{\xi}= \beta (\xi) v_{\xi_1}$ and $v_{\xi_0} = \beta(\xi_0) v_{\xi_1}$ for some non-zero scalars $\beta(\xi), \beta(\xi_0)$. It follows that $\beta(\xi_0)^{-1} = \alpha (\xi_1)$ and $v_{\xi}= \beta (\xi) v_{\xi_1} =  \beta (\xi) \alpha (\xi_1) v_{\xi_0}  =  \beta (\xi) \beta(\xi_0)^{-1} v_{\xi_0}.$ That is, $v_{\xi}$ is a scalar multiple of $v_{\xi_0}$. \smallskip

We have therefore shown the existence of a mapping $\tau: K_1\to \mathbb{C}$ satisfying $v_{\xi} = \tau(\xi) v_{\xi_0}$ for all $\xi \in K_1$. Combining this fact with the conclusion in \eqref{eq definition property of vxi} we arrive at $$\Phi (\xi\otimes \eta) = \widehat{\xi}\otimes \tau(\xi) v_{\xi_0} (\eta) = \overline{\tau(\xi)}  \widehat{\xi}\otimes  v_{\xi_0} (\eta),$$ for all $\xi\in K_1$, $\eta\in H_1$. Defining $u: K_1\to K_1$ by $u(\xi) = \overline{\tau(\xi)}  \widehat{\xi}$, we get a well-defined bijection which inherits the linearity from that of $\Phi$. This concludes the proof of the first statement. \smallskip

Let us briefly comment the second case. \smallskip

\emph{Step 4.} Let us assume now, by Steps 1 and 2, that for each $\xi \in K_1$, there exists a $\widehat{\eta}\in H_2$, depending on $\xi$, such that $\Phi( L_{\xi}) = R_{\widehat{\eta}}$. Consequently, for each $\eta\in H_1$ there exists a unique $u_{\xi} (\eta)\in K_2$ such that \begin{equation}\label{eq new definition property of uxi}  \Phi (\xi\otimes \eta) = u_{\xi}  (\eta) \otimes \widehat{\eta}.
\end{equation} Now the mapping $u_{\xi}: H_1\to K_2$ is a conjugate-linear bijection --essentially because the mapping $(\xi, \eta) \mapsto \xi \otimes \eta$ is sesquilinear. By repeating or adapting the arguments in the first part of the proof, we find two conjugate-linear bijections $u: H_1\to K_2$ and $v: K_1\to H_2$ such that $$\Phi (\xi\otimes \eta ) = u(\eta) \otimes v(\xi), \hbox{ for all } x\in K_1, \eta\in H_1.$$ 
\end{proof}

It should be commented that there is certain margin to consider weaker hypotheses in the above theorem, but the current statement is enough for our purposes. \smallskip

We have already gathered the required machinery to study bijections preserving triple transition probabilities between subsets of minimal tripotents of two type 1 Cartan factors. 

\begin{theorem}\label{t bijections preserving triple transition pseudo-probabilities type 1 Cartan factors} Let $\Phi : \mathcal{U}_{min}(M) \to \mathcal{U}_{min}(N)$ be a bijective transformation preserving triple transition pseudo-probabilities {\rm(}i.e., $TTP(\Phi(v), \Phi(e))=\varphi_{\Phi(e)} (\Phi(v)) = \varphi_{e} (v) = TTP(v,e),$ for all $e,v$ in $\mathcal{U}_{min} (M)${\rm)}, where $M=B(H_1,K_1)$ and $N= B(H_2,K_2)$ are type 1 Cartan factors with dim$(H_j)$, dim$(K_j)\geq2$. Then there exists a {\rm(}unique{\rm)} triple isomorphism $T: M \to N$ whose restriction to $\mathcal{U}_{min} (M)$ is $\Phi$. 
\end{theorem}

\begin{proof} By applying Theorem \ref{t bijections preserving triple transition pseudo-probabilities} we find a bijective linear mapping $$T_0: soc(B(H_1,K_1))\to soc(B(H_2,K_2))$$ whose restriction to $\mathcal{U} (B(H_1,K_1))$ is $\Phi$. Theorem \ref{t extension of OmladicSemrl} asserts that one of the next statements holds: \begin{enumerate}[$(a)$]\item There are bijective linear mappings $u:K_1\to K_2,$ and $v:H_1\to H_2$ such that $\Phi (\xi \otimes \eta) = u(\xi) \otimes v(\eta)$ {\rm(}$\xi\in K_1,\eta \in H_1${\rm)}.
\item There are bijective conjugate-linear mappings $u:H_1\to K_2,$ $v:K_1\to H_2$ such that $\Phi (\xi \otimes \eta)= u(\eta) \otimes v(\xi)$ {\rm(}$\xi\in K_1,$ $\eta \in H_1${\rm)}.
	\end{enumerate}	

The hypothesis affirming that $\Phi$ maps minimal tripotents to minimal tripotents can be now applied to deduce that in any of the previous cases the mappings $u$ and $v$ are isometries. Therefore, $u$ and $v$ are surjective linear isometries in case $(a)$ and surjective conjugate-linear isometries in case $(b)$. So, by defining $T(x) = u x v^*$ and $T(x) = u x^* v^*$ ($x\in B(H_1,K_1)$) in cases $(a)$ and $(b)$, respectively, we get the desired triple isomorphism. 
\end{proof}

\smallskip\smallskip

\textbf{Acknowledgements} Author partially supported by MCIN/AEI/FEDER ``Una manera de hacer Europa'' project no.  PGC2018-093332-B-I00, Junta de Andaluc\'{\i}a grants FQM375, A-FQM-242-UGR18 and PY20$\underline{\ }$00255, and by the IMAG--Mar{\'i}a de Maeztu grant CEX2020-001105-M/AEI/10.13039/501100011033.

\end{document}